\documentclass[12pt]{article}
\date{}
\usepackage{amsmath}
\usepackage{latexsym}
\usepackage{amsfonts}
\usepackage{amssymb}
\usepackage{amscd}
\usepackage{theorem}
\usepackage{pstricks}
\usepackage{pst-node}
\usepackage{latexsym}
\newtheorem{Lemma}{Lemma}
\newtheorem{Proposition}[Lemma]{Proposition}

\newtheorem{Theorem}[Lemma]{Theorem}

\newtheorem{OP}[Lemma]{Open Problem}

\newtheorem{Corollary}[Lemma]{Corollary}         

\newcommand{\Real}{\ensuremath{\mathbb{R}}}

\newcommand{\bck}{\!\!\!}

\newcommand{\GG}{\mbox{${\cal G}$}}         
\newcommand{\FF}{\mbox{${\cal F}$}}         

\newcommand{\eps}{\varepsilon}

\newcommand{\bA}{{\bf A}}         
\newcommand{\bB}{{\bf B}}         
\newcommand{\bE}{{\bf E}}         
\newcommand{\CX}{{\cal X}}         
\newcommand{\bx}{{\bf x}}         
\newcommand{\by}{{\bf y}}         
  
\newcommand{\bu}{{\bf u}}         
\newcommand{\bv}{{\bf v}}         
         
\newcommand{\bee}{{\bf e}}

\newcommand{\rD}{{\rm D}}

\newcommand{\Levy}{{L$\rm \acute{e}$vy}}         
\newcommand{\ERRW}{edge reinforced random walk}

\newtheorem{rema}{Remark}[section]

\newenvironment{proof}{{\sc proof:~}}{\hfill $\Box$  }

\renewcommand{\le}{\leqslant}         
\renewcommand{\ge}{\geqslant}         
         
\newcommand{\bal}{\begin{align*}}         
\newcommand{\eal}{\end{align*}}         
\newcommand{\beq}{\begin{eqnarray*}}         
\newcommand{\eeq}{\end{eqnarray*}}         
\newcommand{\bte}{\begin{theorem}}         
\newcommand{\ete}{\end{theorem}}         
\newcommand{\bl}{\begin{lemma}}         
\newcommand{\el}{\end{lemma}}         
\newcommand{\bde}{\begin{description}}         
\newcommand{\ede}{\end{description}}         
         
\newcommand{\bc}{\begin{cases}}         
\newcommand{\ec}{\end{cases}}         
\newcommand{\bpr}{\begin{proof}}         
\newcommand{\epr}{\end{proof}}         
\newcommand{\bco}{\begin{Corollary}}         
\newcommand{\eco}{\end{Corollary}}

\newcommand{\iy}{\infty}         
\newcommand{\tx}{\text}         
\newcommand{\R}{\ensuremath{\mathbb{R}}}

\newcommand{\Z}{\ensuremath{\mathbb{Z}}}         
\newcommand{\N}{\ensuremath{\mathbb{N}}}

\newcommand{\F}{\mathcal{F}}

\newcommand{\Es}{\mathbb{E}}         
\newcommand{\Pb}{\mathbb{P}}         
\newcommand{\Eb}{\mathbb{E}}         

\newcommand{\al}{\alpha}         
         
\newcommand{\be}{\beta}

\newcommand{\ki}{\chi}         
            
\newcommand{\ka}{\kappa}

\newcounter{saveeqn}%

\newcommand{\reseteqn}{\setcounter{equation}{\value{saveeqn}}%
\renewcommand{\theequation}{\arabic{equation}}}%

\begin{document}         
         
\title{What is the difference between a square and a triangle?}

\author{Vlada Limic and Pierre Tarr\`es}         
\maketitle         
         
\begin{abstract}  
We offer a reader-friendly introduction to the 
attracting edge problem (also known as the ``triangle conjecture")
and its most general current
 solution of Limic and Tarr\`es
(2007).
Little original research is reported; rather
this article ``zooms in'' to describe the essential
characteristics
of two different techniques/approaches 
verifying the almost sure existence
of the attracting edge for the strongly edge reinforced 
random walk (SERRW) on a square.
Both arguments extend straightforwardly to the  
SERRW on even cycles.
Finally, we show that the case where 
 the underlying graph is a triangle cannot be studied by a simple  
modification of either of the two techniques.
\end{abstract}

\vspace{0.2in}         
         
         

{\em Running Head:}         
Difference between a square and a triangle
         
         
\section{Introduction}         
\label{Intro}
We briefly describe the general setting introduced, for example, in \cite{limtar}.
Let $\GG$ be a connected graph          
 with set of vertices $V=V(\GG)$, and set of         
(unoriented)         
edges $E=E(\GG)$.         
The only assumption on the graph is that          
each vertex has at most $D(\GG)$ adjacent vertices         
(edges), for some $D(\GG) < \infty$,         
so that $\GG$ is of bounded degree.         
                  
Call two vertices $v,v'$ {\em adjacent} ($v\sim v'$ in symbols)         
if there exists an edge, denoted by $\{v,v'\}=\{v',v\}$,         
connecting them.         

Let $\N=\{0,1,\ldots\}$, and
let $W: \N \longrightarrow (0,\infty) $  be the {\em reinforcement weight$\,$} function.         
Assume we are          
given initial {\em edge weights} $X_0^e\in\N$ for all $e \in E$,         
such that $\sup_e X_0^e<\infty$.
Let $I_n$ be a $V$-valued random variable, recording the         
position of the particle at time $n\in \N$.   Set $I_0:=v_0$ for some $v_0\in\GG$.       
Let $(\FF_n,\,n\geq 0)$ be the filtration generated by $I$.

The {\ERRW} (ERRW) on $\GG$         
evolves according to a random dynamics with the following         
properties:         
         
(i) if currently at vertex         
$v \in \GG$, in the next step the particle jumps         
to a nearest neighbor of $v$,         
         
(ii)          
the probability of a jump from $v$ to $v'$ at time $n$ is {\em ``$W$-proportional''}         
to the number of previous traversals of the edge connecting $v$ and $v'$, that is,
\[
  \Pb(I_{n+1}=v'|\FF_n)  1_{\{I_n=v\}}= \frac{W(X_n^{\{v,v'\}})}{
  \sum_{w\sim v}W(X_n^{\{v,w\}})} 1_{\{I_n=v\sim v' \}},
\]
where $X_n^e$, $e\in E(\GG)$ equals          
\[
X_n^e=X_0^e + \sum_{k=0}^{n-1} 1_{\{\{I_k,I_{k+1}\}=e \}}.
\]
We recommend a recent survey by Pemantle \cite{pem_survey} 
as an excellent overview of 
processes with reinforcement:  results, techniques, open problems and applications.

Let (H) be the following condition on $W$:         
\newcounter{Heqn}%
\setcounter{Heqn}{0}         
\setcounter{saveeqn}{\value{equation}}         
\begin{alpheqnh}         
\begin{equation}         
\label{H}         
\sum_{k\in\N}\frac{1}{W(k)}<\iy.         
\end{equation}         
\end{alpheqnh}         
\reseteqn
$\bck$We call any {\ERRW} corresponding to $W$ that satisfies (H) a {\em strongly} {\ERRW}.
Denote by 
\[
A:=\{\exists\ n: \{I_k, I_{k+1}\} =\{I_{k+1}, I_{k+2}\},k\geq n\}
\]
 the event that eventually the particle traverses a
single (random) edge of the graph. On $A$ we call that edge
the {\em attracting edge}.    
It is easy to see that          
(H) is         
the necessary and sufficient condition          
for          
\[         
\Pb(\{I_n,I_{n+1}\}=\{I_0,I_1\} \mbox{ for all } n)>0.        
\]         
This implies that (H) is necessary and sufficient          
for $\Pb(A)>0$. The necessity can be seen by splitting $A$ into a countable union of events, where each  corresponds to getting attracted to a specific edge after a particular time with a specific configuration of weights on the neighbouring edges.   
Since $A $ is a tail event, it seems natural to wonder whether
\setcounter{saveeqn}{\value{equation}}         
\begin{alpheqna}         
\begin{equation}         
\label{A}         
\Pb(A)=1         
\end{equation}         
\end{alpheqna}         
\reseteqn
$\bck$holds.          
The authors studied this problem in \cite{limtar}, 
and concluded that, under  additional technical assumptions,
(\ref{H}) implies (\ref{A}).
In particular, 
\begin{Theorem} [\cite{limtar}, Corollary 3]
\label{limtar}
If $\GG$ has bounded degree and if
$W$ is non-decreasing, then (\ref{H}) implies (\ref{A}).
\end{Theorem}

We denote by $\GG_{l}$ the {\em cycle } of length $l$, 
with vertices $\{0,1,\ldots,l-1\}$
and edges
\[
e_i = \{i,i+1\}, 
 \ i =0,\ldots,l-1,
\]
where $l\geq 3$, and where the addition is done modulo $l$. 

 Let us now  concentrate on the case where the
underlying graph $\GG $ is the square $\GG_4$.
The next two sections demonstrate two different techniques
of proving the following claim.
\begin{Theorem}  
\label{Tsquare}
If $\GG =\GG_4 $, then (\ref{H}) implies (\ref {A}).
\end{Theorem}
In fact we will concentrate on a somewhat simpler claim whose proof can be 
``recycled'' (as we henceforth discuss) 
in order to arrive to the full statement of Theorem \ref{Tsquare}.
\begin{Proposition}  
\label{Psquare}
If $\GG =\GG_4 $, then (\ref{H}) implies 
\[
\Pb(\tx{all four edges are traversed infinitely often})=0.
\]
\end{Proposition}

In Section \ref{S:explain} we discuss the reasons why
these techniques which are well-suited for
$\GG =\GG_4$, or any graph of bounded degree without an odd cycle, cf.~\cite{sel94}
or \cite{lim03},
do not extend to the setting where
$\GG =\GG_3$ is a triangle.
In fact, the following "triangle conjecture" is still open 
in its full generality
(cf.~Theorem \ref{Tsquare} and Theorem \ref{limtar})
where $W$ is a general (irregular) weight function satisfying (\ref{H}): 
\begin{OP} 
If $\GG =\GG_3 $, then (\ref{H}) implies (\ref {A}).
\end{OP}

\section{A continuous time-lines technique }
\label{Stime}
This technique adapts a construction due to Rubin, and was invented by Davis \cite{dav90} and 
Sellke \cite{sel94}. It played a key role in his proof
of the attracting edge property on $\Z^d $, and was also used by
Limic \cite{lim03} 
in order to simplify the attracting edge problem on graphs of bounded degree to the same problem
on odd cycles.
Denote for simplicity
\[
e_i := \{i,i+1\}, 
 \ i =0,\ldots,3,
\]
where addition is done modulo $4$.
For each $i=0,\ldots,3$ and $k\geq 1$ let $E_k^i$ be an exponential random variable with mean $1/W(k)$,
such that $\{E_k^i,i=0,\ldots,3,k\geq 1\}$ is a family of independent random variables.
Denote by
\[
T_n^i:=\sum_{k=X_0^{e_i}}^{X_0^{e_i}+n} E_k^i, \ n\geq 0,
\
T_\infty^i:=\sum_{k=X_0^{e_i}}^\infty E_k^i, \ i=0,\ldots,3.
\]
Note that the random variables $T_n^i$,  $T_\infty^i$, $i=0,\ldots,3$, $n\in\N$, are continuous, independent and finite almost surely
(the last property is due to assumption (H)).
In Figure 1, the $T_n^i$ are shown as dots, and 
the ``limits'' $T_\infty^i$, $i=0,\ldots,3$ are indicated.

\psset{xunit= 1.05cm,yunit= 1cm}
\pspicture(0,-2.5)(15,2)
\psset{linewidth=.5pt}
\psline{|->}(0,-2)(8.5,-2)
\psline{|->}(0,-1)(8.5,-1)
\psline{|->}(0,0)(8.5,0)
\psline{|->}(0,1)(8.5,1)
\pscircle*(2,-2){.05} \pscircle*(3.5,-2){.05} \pscircle*(4.8,-2){.05}
\pscircle*(5.8,-2){.05}\pscircle*(6.2,-2){.05}\pscircle*(6.5,-2){.05}
\pscircle*(6.7,-2){.05}
\pscircle*(1,-1){.05} \pscircle*(2.3,-1){.05} \pscircle*(3.32,-1){.05}
\pscircle*(4.1,-1){.05}\pscircle*(4.3,-1){.05}\pscircle*(4.4,-1){.05}
\pscircle*(4.5,-1){.05}
\pscircle*(1.1,0) {.05} \pscircle*(2.5,0) {.05} \pscircle*(3.7,0){.05}
\pscircle*(4.45,0){.05}\pscircle*(4.8,0){.05}
\pscircle*(4.91,0){.05}
\pscircle*(0.9,1){.05} \pscircle*(1.6,1){.05} \pscircle*(3.1,1){.05}
\pscircle*(4,1){.05}\pscircle*(4.15,1){.05}
\rput[u](0,-2.4){$0$}
\rput[u](4.8,0.6){$T_\infty^0$}
\rput[u](5.6,-.4){$T_\infty^1$}
\rput[u](5.2,-1.4){$T_\infty^2$}
\rput[u](7.4,-2.4){$T_\infty^3$}
\rput[u](1,0.6){$T_0^0$}
\rput[u](1.3,-.4){$T_0^1$}
\rput[u](1.2,-1.4){$T_0^2$}
\rput[u](2,-2.4){$T_0^3$}
\rput[u](1.7,0.6){$T_1^0$}
\rput[u](2.5,-.4){$T_1^1$}
\rput[u](2.3,-1.4){$T_1^2$}
\rput[u](3.5,-2.4){$T_1^3$}
\rput[b](4.5,1){$\cdots$}
\rput[b](5.3,0){$\cdots$}
\rput[b](4.9,-1){$\cdots$}
\rput[b](7.1,-2){$\cdots$}
\rput[u](9.5,1.5){time-line of}
\rput[l](9.35,1){$e_0$}
\rput[l](9.35,0){$e_1$}
\rput[l](9.35,-1){$e_2$}
\rput[l](9.35,-2){$e_3$}
\endpspicture

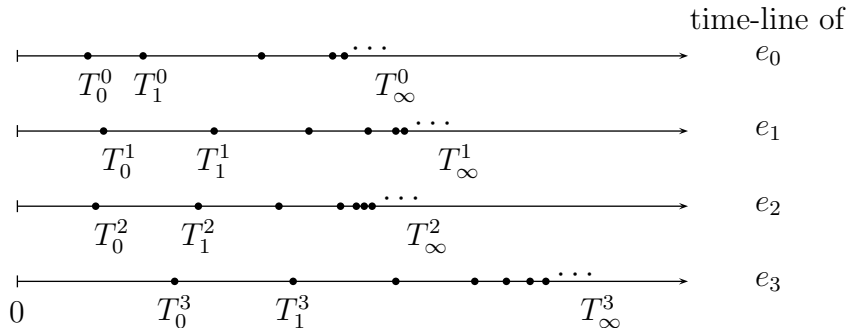
\centerline{Figure 1}

\vspace{0.2cm}
\noindent

Here is how one can construct a realization of the {\ERRW} on $\GG_4$ from the above data,
or (informally) from the figure.
Given the current position of the walk, simultaneously erase 
(at rate $1$) the two time-lines of the incident edges in the chronological direction 
until encountering the next dot belonging to either of the time-lines.
At this point, the walk steps into a new vertex by 
traversing the edge 
that corresponds to 
the time-line 
containing the
dot.
The procedure continues inductively.  

We next explain  why this construction indeed 
leads to a realization of 
the {\ERRW}, by considering carefully the first three steps,.
Assume for concreteness that the initial position is vertex $0$ incident to 
the edges $e_0$ and $e_3$.
The time-lines of $e_0$ and $e_3$ 
are erased until 
the minimum of $T_0^0=E_{X_0^{e_0}}^{0}$ and $T_0^3=E_{X_0^{e_3}}^{3}$.
In the figure this minimum happens to be $T_0^0$.
Thus the particle moves from $0$ to $1$ (traversing edge $e_0$) in the first step.
Due to the properties of exponentials, the probability of this move is exactly
$W(X_0^{e_0})/(W(X_0^{e_0})+W(X_0^{e_3}))$.
The two incident edges to the current position $I_1$ are now $e_0$ and $e_1$.
Continue by simultaneous erasing 
(the previously non-erased parts of) time-lines corresponding to $e_0$
and $e_1$ until the next dot.
In the figure, the dot again appears on the line of $e_0$. 
Hence the particle  traverses the edge $e_0$ in the second step and therefore jumps back to 
vertex $0$. Note that again the probability of this transition matches the one of the {\ERRW}.
Continue by the 
simultaneous erasure 
of time-lines corresponding to $e_0$ and $e_3$.
Based on the figure, the particle makes the third step across the edge $e_3$, since the 
(residual) length of the interval on the time-line of $e_3$ until 
$T_0^3$ is smaller than $T_2^0- T_1^0=E_{X_0^{e_0}+2}^1$.
The 
memoryless property of the exponential
distribution insures that the (residual) length of the interval until 
$T_0^3$ is again distributed as exponential (rate $W(X_0^{e_3})$) random variable,
 independent of all other data.
 Hence, the transition probability again matches that of the ERRW.

Note that the above construction can be done with any number
$l\geq 3$ of time-lines (corresponding to the length $l$ of the 
underlying circular graph), 
and we make use of this generalization in 
Section \ref{S:explain}.

As a by-product of the above construction, a continuous-time version  
of the {\ERRW} emerges, where the particle makes the  jumps exactly at times when the dots are encountered.
More precisely,
if we denote by  $\widetilde{I}(s)$ the position of the particle at time $s$ 
and if $\tau_0=0$ and $0<\tau_1<\tau_2<\ldots$ are the successive jump times 
of $\widetilde{I}$,
then the 
(discrete time) ERRW constructed above,
and the continuous-time version are coupled so that
\[
I_k \equiv \widetilde{I}(\tau_k), \ k \geq 0.
\]

It is worth noting that this continuous-time version is analogous to the Harris construction of a continuous-time Markov chain from the discrete one, yet it is different since the parameters of the exponential clocks vary. In particular, under assumption (\ref{H}), the total time of evolution of the continuous-time random walk is finite. 

Consider the total time of evolution for the continuous time walk,
\[
T := \lim_{k\to \infty} \tau_k.
\]
Note that at any time $s\geq 0$ the particle is incident to one of the edges $e_0$ and $e_2$, and
equally it is  incident to one of the edges $e_1$ and $e_3$,
hence
\begin{eqnarray}
T
&=&\sum_{i=0, \,i \tx{ even}}^3 \mbox{
total 
time spent on boundary vertices of }e_i 
\nonumber\\
&=&\sum_{i=0, \,i \tx{ odd}}^3\mbox{total 
time spent on boundary vertices of }e_i.
\nonumber
\end{eqnarray}
Note that
\begin{eqnarray}
\label{Eeventsame}
& &\{
\tx{all four edges are traversed infinitely often}
\}\\
& &=
\nonumber
\{\mbox{the time-line of $e_i$ is erased up to time $T_\infty^i$ for each $i=0,\ldots,3$}\}
\\
& &\subset\label{Eequal}
\{T= T_\infty^0 + T_\infty^2 = T_\infty^1 + T_\infty^3
\}.
\end{eqnarray}
However, due to the independence and continuity of $T_\infty^0 + T_\infty^2$ 
and $T_\infty^1 + T_\infty^3$, the identity (\ref{Eequal})
happens with probability $0$.
We conclude that (\ref{Eeventsame}) happens with probability $0$, and therefore that
Proposition \ref{Psquare} holds.

In order to obtain the proof of Theorem \ref{Tsquare} now note that there are 
essentially three possibilities
remaining for the asymptotic evolution: the {\ERRW} visits infinitely often 
either one, or two adjacent, or three edges.
In the latter two cases, 
there is at least one vertex $j$ such that both edges
$e_{j-1}$ and $e_j$ are traversed infinitely often.
Moreover, after finitely many steps, 
every excursion from $j$ 
starts and ends
with the same 
edge.
Now one can measure the time spent at site $j$ 
from two perspectives: that of waiting to traverse edge $e_{j-1}$, and 
that of waiting to traverse edge $e_j$.
The reader will quickly construct a variation 
to the above argument (alternatively, consult \cite{sel94} or \cite{lim03}) 
determining that a branching vertex exists with probability $0$.

Note that this continuous time-lines technique still works on even cycles 
$\GG_{2k}$.
Indeed, given the continuous-time realization of the {\ERRW} constructed above, we observe that on the event that all edges 
are visited infinitely often, 
\begin{equation}
\label{eqtimes}
T=\sum_{i=0, i\tx{ even}}^{2k-1}T_\iy^i=\sum_{i=0, i\tx{ odd}}^{2k-1}T_\iy^i,
\end{equation}
where
$T:= \lim_{k\to \infty} \tau_k$ is the total time of evolution for 
the walk.
As before, (\ref{eqtimes}) is  
a consequence of the fact that
$T$ equals the total time spent 
on both the boundary of even and the boundary of odd edges.
Now, due to independence and continuity of 
$\sum_{i=0, i\tx{ even}}^{2k-1}T_\iy^i$ and $\sum_{i=0, i\tx{
    odd}}^{2k-1}T_\iy^i$,  the identity (\ref{eqtimes}) happens with
probability $0$ so that, almost surely, at least one of the edges in
the cycle is  visited
only finitely often, and we conclude (\ref{A}) as in the case of the square. 
\section{A martingale technique}
\label{Smart}
Let, for all $n\in\N$, 
$$W^*(n):=\sum_{k=0}^{n-1}\frac {1}{W(k)},$$
with the convention that $W^*(0):=0$.

Assume the setting of Proposition \ref{Psquare}.
For all $n\in\N$, $i=0,\ldots,3$, define the processes
\begin{align}
\label{yn}
Y_n^\pm(i)&:=\sum_{k=0}^{n-1} \frac{1_{\{I_k= i,I_{k+1}=i\pm1\}}}{W(X_k^{\{i,i\pm1\}})} 
\\
\label{kan}
\kappa_n^i &:= Y_n^+(i)-Y_n^-(i)
\end{align}

Clearly, $\kappa_\cdot^i$ is measurable with respect to the filtration 
$(\FF_n,\,n\geq 0)$. 
Moreover, it is easy to check that $(\kappa_n^i,\ n\geq 0)$ is a martingale : 
on $\{I_n=i\}$, 
$E(\kappa_{n+1}^i-\kappa_n^i|\FF_n)$ is equal to
\[
\frac{1}{W(X_n^{e_i})} \frac{W(X_n^{e_i})}{W(X_n^{e_i})+W(X_n^{e_{i-1}})}
-
\frac{1}{W(X_n^{e_{i-1}})} \frac{W(X_n^{e_{i-1}})}{W(X_n^{e_i})+W(X_n^{e_{i-1}})}
=0.
\]
Therefore the process
\begin{equation}
\label{Ekappa}
\kappa_n := \kappa_n^0 - \kappa_n^1 + \kappa_n^2 -\kappa_n^3
+\sum_{i=0}^3(-1)^i W^*(X_0^{e_i}),
\end{equation}
is also a martingale.
Due to assumption (H), each of the four processes $\kappa_.^i$ 
is a difference of bounded non-decreasing processes, and
therefore has an almost sure limit as $n\to \infty$.
Hence denote by $\kappa_\infty$ the finite limit $\lim_{n\to \infty} \kappa_n$.

Now
\begin{equation}
\label{eq:wstar}
\ka_n=\sum_{i=0}^3(-1)^i W^*(X_n^{e_i}).
\end{equation}
This implies that 
\[
\{\tx{all four edges are traversed infinitely often}\}
\subset
\{\kappa_\infty =0\},
\]
so that it suffices to show 
\begin{equation}
\label{Eeventconsid}
\Pb({\cal A}_\infty)=0,
\end{equation}
where \[
{\cal A}_\infty:=\{\tx{all four edges are traversed infinitely often}\} \cap \{ \kappa_\infty=0 \}.
\]

In order to prove (\ref{Eeventconsid}), we now analyze carefully the
variance of the increments of the martingale $(\ka_n)_{n\in\N}$ 
(decreasing to $0$, due to (H)), which will enable us to prove the
nonconvergence of this martingale to $0$ a.s.~on the event that all
edges are  visited infinitely often. This technique adapts an argument proving almost sure nonconvergence  towards unstable points of stochastic approximation algorithms, introduced by Pemantle \cite{pem90} and generalized by Tarr\`es \cite{tarres1,tarres2}.

Fix large $n$, and note that
\begin{eqnarray}
\nonumber & &\Eb((\kappa_{n+1})^2-(\kappa_n)^2|\FF_n)=\Eb((\kappa_{n+1}-\kappa_n)^2|\FF_n) \\
& &\label{Eleveln}
= \Eb\!\!\left(\left.\sum_{i=0}^3\frac{1_{\{I_{n}=i,I_{n+1}=i+1\}}}{(W(X_n^{e_i}))^2}
+ \frac{1_{\{I_{n}=i,I_{n+1}=i-1\}}}{(W(X_n^{e_{i-1}}))^2}\right|\FF_n\!\!\right).
\end{eqnarray}
From now on abbreviate
\[
\alpha_n := \sum_{j=X_n^*}^\infty \frac{1}{(W(j))^2},
\]
where $X_n^* = \min_{i=0,\ldots,3} X_n^{e_i}$.
For $\eps>0$, define the stopping time
\begin{equation}
\label{ES}         
S:=\inf\{k\ge n:\,|\ka_k|>\eps\sqrt{\al_n}\}.
\end{equation}
Since 
\[
(\kappa_S)^2 -(\kappa_n)^2=\sum_{k=n}^\infty ((\kappa_{k+1})^2-(\kappa_k)^2) 1_{\{ S>k\}},\]
by nested conditioning we obtain
\[
\Eb((\kappa_S)^2 -(\kappa_n)^2|\FF_n)=\Eb (\sum_{k=n}^\infty \Eb[(\kappa_{k+1})^2-(\kappa_k)^2|\FF_k]
 1_{\{ S>k\}}
|\FF_n ),
\]
so that, due to (\ref{Eleveln}), we obtain
\begin{eqnarray}
\Eb((\kappa_S)^2 \bck &-&\bck(\kappa_n)^2 |\FF_n)= 
\Eb\left[
\sum_{k=n}^{S-1}\sum_{i=0}^3\frac{1_{\{I_{k}=i,I_{k+1}=i+1\}}}{(W(X_k^{e_i}))^2}
+ \frac{1_{\{I_{k}=i,I_{k+1}=i-1\}}}{(W(X_k^{e_{i-1}}))^2}
\big|\FF_n\right]\nonumber\\
&=&
\label{EsqlowerA}
\sum_{i=0}^3
\Eb\left[
\sum_{k=X_n^{e_i}}^{X_{S}^{e_i}-1}
\frac{1}{W(k)^2}
\big|\FF_n\right]\geq \alpha_n \Pb({\cal A}_\infty \cap \{S=\infty\}|\FF_n) .
\end{eqnarray}
However, $\kappa_S=0$ on $\{S=\infty\}\cap\cal A_\iy$, also $|\ka_S|=|\ka_\iy|\le\eps\sqrt{\al_n}$ on $\{S=\iy\}$ and, on $\{S<\iy\}$,  $|\kappa_S| \leq (1+\eps)\sqrt{  \al_n}$ since the over(under)shoot of  
$\kappa$ at time $S$ is bounded by a term of the
type $1/W(l)$ for some random $l\geq X_n^*$, 
so in particular it is bounded by $\sqrt{\alpha_n}$. Hence 
\begin{equation}
\label{Elevelb}
\Eb((\kappa_S)^2 |\FF_n)\leq
\Eb((\kappa_S)^21_{\{S<\infty\}\cup{\cal A}_\iy^c} |\FF_n) \leq  
(1+\eps)^2\al_n\Pb(\{S<\infty\}\cup{\cal A}_\iy^c |\FF_n).
\end{equation}

Letting $p:=\Pb({\cal A}_\infty \cap \{S=\infty\}|\FF_n)$, we conclude by combining inequalities 
(\ref{EsqlowerA})  and (\ref{Elevelb}) that  
$p  \leq 
(1+\eps)^2 (1-p),$
or equivalently 
\begin{equation}
\label{Ealmostdone}
\Pb({\cal A}_\infty \cap \{S=\infty\}|\FF_n) =p\leq (1+\eps)^2/(1+(1+\eps)^2)<1,
\end{equation}
almost surely.

It will be convenient to let  $\eps=5$.
Then note that the shifted process 
$(\kappa_{S+k},k\geq 0)$ is again a martingale with respect to 
the filtration $\tilde{\FF}_k := \FF_{S+k}$.
Moreover, due to (\ref{Eleveln}), we have that
\[
\Eb((\kappa_\infty-\kappa_S)^2|\FF_S)\leq 4\alpha_S \leq 4\alpha_n,
\]
so that by the Markov inequality, a.s.~on $\{S<\infty\}$,
\[
\Pb({\cal A}_\infty|\FF_S)\leq \Pb(|\kappa_\iy-\kappa_S|> 5 \sqrt{\alpha_n}|\FF_S)\leq \frac{4 \alpha_n}{25 \alpha_n}=
\frac{4}{25},
\]
thus
\[
\Pb({\cal A}_\infty^c|\FF_n)\geq \Eb [\Pb({\cal A}_\infty^c|\FF_S)1_{\{S<\infty\}}|\FF_n]
\geq \frac{21}{25} \Pb(S<\infty|\FF_n).
\]
Note that (\ref{Ealmostdone}) now implies 
\[
\Pb({\cal A}_\infty^c|\FF_n) \left(1+ \frac{25}{21}\right) \geq 
\Pb({\cal A}_\infty^c|\FF_n) + \Pb(S<\infty|\FF_n) \geq
1-(1+\eps)^2/(1+(1+\eps)^2),
\]
so finally
\[
\Pb({\cal A}_\infty^c|\FF_n)\geq c,\]
almost surely for some constant $c>0$.
By the {\Levy} 0-1 law, we conclude that Proposition \ref{Psquare} holds. 

In order to prove Theorem \ref{Tsquare} we can proceed as in the previous section
to show that no branching point is possible. 
In particular, we consider $i,j\in\{0,\ldots,3\}$ such that $j\neq i,i-1$,
and events of the form 
\[
\{e_i \tx{ and } e_{i-1} \tx{ both traversed i.o.}\}\cap \{e_j \tx{ not visited after time } n\},
\]
for some finite $n$, and then use an appropriate modification of
$(\kappa_k^i,\, k\geq n)$ that would have to converge to a particular limit 
on the above event, and show in turn that this convergence occurs with 
probability $0$.

Note that again this martingale technique extends in the more general setting of even cycles $\GG_{2k}$. Indeed, 
let $Y_n^\pm(i)$ and $\ka_n^i$ be defined as in (\ref{yn}) and
(\ref{kan}) and, let 
$$
\kappa_n := \sum_{i=0}^{2k-1}(-1)^i\ka_n^i+\sum_{i=0}^{2k-1}(-1)^i W^*(X_0^{e_i}).$$

As in equation (\ref{eq:wstar}), 
$$\ka_n=\sum_{i=0}^{2k-1}(-1)^i W^*(X_n^{e_i}),$$
so that 
\[
\{\tx{all edges are traversed infinitely often}\}
\subset
\{\kappa_\infty =0\}.
\]
The study of the variances of the martingale increments explained in
Section \ref{Smart} yields similarly that $\Pb(\{\ka_\iy=0\})=0$. 
Hence, almost surely, at least one of the edges in the cycle is 
visited only finitely often and, as before,
an 
adaptation of this argument implies (\ref{A}). 

\section{Comparing square and triangle}
\label{S:explain}
In order to additionally motivate our interest in the evolution of {\ERRW} on cycles, we recall that the
continuous 
time-line technique can be adapted in order to prove that, on any
graph of bounded degree, almost surely, the strongly {\ERRW}
either satisfies (\ref{A}) or it eventually keeps traversing
infinitely often all the edges of a (random) odd sub-cycle.
The argument was given by Limic in \cite{lim03}, Section 2, using  
graph-based techniques. The martingale method  could be used in a 
similar way to prove the above fact.
In view of this, note that solving the attracting edge problem on odd 
cycles is necessary and sufficient for obtaining the solution on
general bounded degree graphs.

The aim of this section is to explain why the continuous time-line and martingale techniques 
do not extend easily to the setting where $\GG$ is an odd cycle
(e.g., a triangle). 

\subsection{Odd versus even in the time-line technique}
The argument in the setting of even cycles relied on the existence of the  non-trivial linear identity (\ref{eqtimes}) involving independent continuous random variables. We are going to argue next that no such non-trivial linear relation (and in fact no non-linear smooth relation either) can hold with positive probability in the odd case. 

Fix $l \geq 3$ and consider the set 
$$
{\cal X}:=\{
(x_k^i)_{k\in \N, i\in\{0,\ldots,l-1\}}: 
 \forall i\in \{0,\ldots,l-1\}, \forall k\geq 0, x_k^i>0 , \sum_m x_m^i < \iy\}.
$$
Given $\bx= (x_k^i)_{k\in \N, i\in\{0,\ldots,l-1\}} \in\CX$, define
$$
t_n^i\equiv t_n^i(\bx):=\sum_{k=0}^n x_{k}^i, n\geq 0, \ t_\iy^i\equiv t_\iy^i(\bx):=\sum_{k=0}^\iy x_k^i,\ i=0,\ldots,l-1.
$$
After placing dots at points $\,t_0^i< t_1^i< \ldots $ 
on the time-line of $e_i$, $i=0,\ldots,l-1$, 
and fixing the starting position $\iota_0$ 
one can perform, as in Section \ref{Stime}, the 
time-line construction of the (now deterministic) walk, driven by $\bx$, evolving in
continuous time. 
If at any point the erasing procedure encounters more than one dot
 (on two or more different time-lines) simultaneously,
choose to step over the edge corresponding to one of these time-lines
in some prescribed way, for example, to the one having the smallest index.
Denote by 
  $s^{i}:=s^{i}(\bx,\iota_0)$ the total time 
this deterministic walk spends visiting 
vertex $i$.
Similarly, denote by 
\begin{equation}
\label{Elocrel}
t^{e_i}=t^{e_i}(\bx,\iota_0):= s^{i}(\bx,\iota_0)+ s^{i+1}(\bx,\iota_0)
\end{equation} 
the total time that this deterministic walk spends waiting on the boundary 
vertices $i,i+1$ of $e_i$.
Of course, $t^{e_i}(\bx,\iota_0)\leq t_\iy^i(\bx)$, where the equality holds 
if and only if $e_i$ is traversed infinitely
often.
In the case of even $l$ the identity
\[
\sum_{j=0,j \tx{ even}}^{l-1}  t^{e_j}(\bx,\iota_0)= \sum_{j=0,j \tx{ odd}}^{l-1}  t^{e_j}(\bx,\iota_0), \ \bx \in \CX, 
\] 
lead to (\ref{eqtimes}) and was the key for showing that (\ref{A}) occurs. 
Let 
$$y\equiv y(\bx,\iota_0):=\left(\begin{array}{c}
s^0\\
\vdots\\
s^{l-1}\
\end{array}
\right), ~
z\equiv z(\bx,\iota_0):=\left(
\begin{array}{c}
t^{e_0}\\
\vdots \\
t^{e_{l-1}}
\end{array}
\right), 
$$
and 
$$
M^{(l)}:=(\ki_{\{i,i+1\}}(j))_{0\le i,j\le l-1},
$$
where $\ki_B$ denotes a characteristic function of a set $B$, and 
the addition is done modulo $l$, for instance 
$$M^{(5)}=\left[
\begin{array}{ccccc}
1&1&0&0&0\\
0&1&1&0&0\\
0&0&1&1&0\\
0&0&0&1&1\\
1&0&0&0&1
\end{array}
\right].
$$

Then (\ref{Elocrel}) states that 
$$z(\bx,\iota_0)=M^{(l)}y (\bx,\iota_0), \ \bx\in \CX.$$
Note that the determinant $\det (M^{(l)})=1-(-1)^l$ can easily be computed explicitly since 
$M^{(l)}$ is a circular matrix. Hence
$M^{(l)}$ is a regular matrix if and only if $l$ is odd.
Therefore, for odd $l$ and fixed $\iota_0$, a nontrivial identity 
\begin{equation}
\label{Ezbetac} 
\be \cdot z(\bx,\iota_0)=c,\ \bx\in \CX, \ 
\end{equation} 
for some $\be \in \Real^l\setminus\{0\},\ c\in \Real$, holds if and only if, 
\begin{equation}
\label{Eybetac} 
\be' \cdot \by(\bx,\iota_0) =c,\ \bx\in \CX, 
\end{equation} 
where
$\be'=(M^{(l)})^\tau \be \in \Real^l$ is again $\neq 0$.

Now (\ref{Eybetac}) cannot hold identically on $\CX$, and 
we are about to show a somewhat stronger statement.
Let $\bx\in \CX$ and fix some $r\in (0,\iy)$.
Then for $j\in \{0,\ldots,l-1\}$, let
 $  \eta_r^j (\bx) \equiv\tilde{\bx}^{(j)}:= 
({x}_k^{i,(j)})_{k\in \N, i=0,\ldots,l-1} \in \CX$ be defined as
follows:
if $k\geq 0$, $i\in\{0,\ldots,l-1\}$,
$$
\tilde{x}_k^{i,(\iota_0)}:=
x_k^i+r \ki_{\{(\iota_0,0),(\iota_0-1,0)\}} ((i,k)),
$$
while for $j\neq \iota_0$, if the walk driven by $\bx$ visits site $j$ 
for the first time
by traversing $e_{j-1}$,  let  
$$
\tilde{x}_k^{i,(j)}:=
x_k^i+r \ki_{\{(j,0),(j-1,1)\}} ((i,k)),
$$ 
 otherwise let   
\begin{equation} 
\tilde{x}_k^{i,(j)}:=
x_k^i+r \ki_{\{(j-1,0),(j,1)\}} ((i,k)).
\label{Etreca}
\end{equation}
Note that (\ref{Etreca})
comprises also the case where the walk 
driven by $\bx$ never visits site $j$.

Now we will modify the {\ERRW}  by delaying the first jump out of a 
particular site $j$ by some positive amount $r$, without changing anything
else in the behaviour. Informally, the law of the modified version will be absolutely continuous with respect to the law of the original, and this will lead to a contradiction.  

More precisely, for each fixed $j$, 
consider the two (deterministic time-continuous) walks:
the original one that is driven by $\bx$, and the new one 
that is driven by the
transformed family 
$\eta_r^j(\bx)$. 
It is easy to check that
either neither of the walks visits site $j$, or  
they both do. In the latter case, if we denote respectively by $a(\bx)$ and 
$a(\eta_r^j(\bx))$ the amount of time they spend at site $j$  
before leaving,
then $a(\eta_r^j(\bx))= a(\bx) +r.$  
Everything else in the evolution of the two walks is the
same.
In particular, if the walk driven by $\bx$ ever visits $j$, then
\begin{equation}
\label{Esimpo}
s^j(\eta_r^j(\bx),\iota_0)= s^j(\bx,\iota_0)+r,
\tx{ and } 
s^i(\eta_r^j(\bx),\iota_0)= s^i(\bx,\iota_0),\, i\neq j.
\end{equation}
Now one can simply see that if  
the walk driven by $\bx$ visits all sites at least once,
then any identity (\ref{Eybetac})
breaks 
in any open neighborhood of $\bx$, due to  points $\eta_r^j(\bx)$  contained in it for sufficiently small positive $r$.  

Recall the setting of Section \ref{Stime},  
and to simplify the notation, assume $X_0^{e_i}=0$, $i=0,\ldots,l-1$, and 
specify the initial position $\widetilde{I}_0=v_0 \in \{0,\ldots,l-1\}$.  
The point of the above discussion is that the 
random walk $\widetilde{I}$ is then the walk driven by
a $\CX$-valued random family
$\bE=(E_k^i)_{k\in\N,i=0,\ldots,l-1}$, where the random variables $E_k^i$, $k\geq 0$
are independent exponentials, as specified in Section \ref{Stime}.  
If 
$$
A_{all}:=\{\widetilde{I} \tx{ visits all vertices at least once}\},
$$
then
$
A_{all}= \{ {\bf E} \in {\bA}_{all}\}
$,
where ${\bA}_{all}$ contains all $\bx\in \CX$ such that the 
deterministic walk
driven by $\bx$ visits all vertices at least once.
It is natural to ask whether one or more (up to countably many, note
that this would still be useful) identities of the form
(\ref{Ezbetac}) hold on 
${\bA}_{all}$, with positive probability.
For $l$ even, we know that the answer is affirmative.
For $l$ odd, this is equivalent to asking whether one or 
more (up to countably many)
identities of the form
(\ref{Eybetac}) hold on ${\bA}_{all}$, with positive probability.
Assume 
\begin{equation}
\label{Eposiproba}
\Pb(A_{all} \cap \{\be' \cdot y(\bE,v_0)=c
\})=p(\be',c)>0,
\end{equation}
for $\be',c$ as in (\ref{Eybetac}).
Take $j$ such that $\be'_j \neq 0$ 
(at least one such $j$ exists).
We will assume that $j\neq v_0$, the argument is somewhat similar and simpler 
otherwise.
Denote by ${\bA}_{1}^{j,-}\subset \CX$ the set of all
$\bx$ such that the walk
driven by $\bx$ visits $j$ for the first time by traversing $e_{j-1}$, and let
$A_{1}^{j,-}=\{\bE \in {\bA}_{1}^{j,-}\}$.
In view of (\ref{Eposiproba}), without loss of generality, we may assume  
\begin{equation}
\label{Eposiprobac}
\Pb(A_{all}\cap A_{1}^{j,-} \cap \{\be' \cdot y(\bE,v_0)=c
\})\geq p(\be',c)/2,
\end{equation}
As a consequence of the earlier discussion in the deterministic setting we have
$
A_{1}^{j,-} =\{\bE\in {\bA}_{1}^{j,-}\} \subset \{\eta_r^j(\bE)\in {\bA}_{1}^{j,-}\},
$
and 
$$
A_{all}\cap \{\be' \cdot y(\bE,v_0)=c\}
\subset
\{ \eta_r^j(\bE)\in \bA_{all} \}\cap 
\{\be' \cdot y(\eta_r^j(\bE),v_0)=c+r\be_j'\},
$$
almost surely.
Therefore, 
\begin{equation}
\label{Eposiprobab}
\Pb( \{ \eta_r^j(\bE) \in \bA_{all} \cap {\bA}_{1}^{j,-}\}
 \cap \{\be' \cdot y(\eta_r^j(\bE),v_0)=c+r\be_j'
\})\geq p(\be',c)/2.
\end{equation}
However, $E_k^i$ are continuous and independent random variables,
each
taking values in any interval $(a,b)\subset (0,\iy)$, $a<b\leq \infty$ 
with positive probability.
Moreover, since $E_0^j$ and $E_1^{j-1}$
are exponential (rate $W(0)$ and $W(1)$, respectively),
 one can easily verify that for any cylinder set $\bB\subset \CX$,
\begin{eqnarray}
\nonumber
 \Pb( \eta_r^j(\bE) \in \bB \cap {\bA}_{1}^{j,-})&=&\\
\nonumber
\Pb( \{E_i^k+r \ki_{\{(j,1),(j-1,0)\}}((i,k)),& &\bck\bck\bck\bck  k\geq 0,j=0,\ldots,l-1 \} 
\in \bB \cap {\bA}_{1}^{j,-})\\
&\leq& 
\label{Ecr}
e^{r(W(0)+W(1))} \Pb(\bE \in \bB \cap {\bA}_{1}^{j,-}).
\end{eqnarray}
Now (\ref{Ecr}) and (\ref{Eposiprobab}) imply
\begin{equation}
\label{Elastineq}
\Pb( \{ \bE \in \bA_{all} \}
 \cap \{\be' \cdot y(\bE,v_0)=c+r\be_j'
\})\geq e^{-r(W(0)+W(1))} p(\be',c)/2,
\end{equation}
and this, together with (\ref{Eposiproba}), leads to a contradiction,
since adding (\ref{Elastineq}) over all rational $r\in (0,1)$ would imply
$\Pb(A_{all})=\Pb( \bE \in \bA_{all})=\iy$.

In the absence of a convenient linear identity (\ref{Ezbetac}), the reader might be tempted to look
for non-linear ones.
Yet, the last argument can be extended to a more generalized setting where 
(\ref{Eybetac}) is replaced by
\begin{equation}
\label{Enonlina}
\by(\bx,\iota_0) \in M, \ \bx\in \CX,
\end{equation}
for some $l-1$-dimensional differentiable manifold $M\subset \R^l$.
In particular, this includes the case where 
$F(\by(\bx,\iota_0))=0, \ \bx\in \CX$, for some smooth function $F$ with non-trivial gradient
(see, for example, \cite{spivak} Theorem 5-1).
Indeed, assume that, in analogy to (\ref{Eposiproba}), 
\begin{equation}
\label{EMposiproba}
\Pb(A_{all} \cap \{ \by(\bE,v_0)\in M\})>0.
\end{equation}
Then, since $\by(\bE,v_0)$ is a finite random vector, due to
 the definition of differential manifolds (cf.~\cite{spivak} p.~109),
there exists a point 
$x\in M$, two bounded open sets $U\ni x,V \subset\R^l$, 
and a diffeomorphism $h:U\to V$  such that
\begin{equation}
\label{EmanifU}
\Pb(A_{all} \cap\{\by(\bE,v_0) \in M \cap U\})=: p(M,U)>0,
\end{equation}
and 
\[
h(U\cap M) = V \cap (\R^{l-1}\times \{0\})=\{\bv \in V: v_l=0\}.
\]
Denote by $\bee_j$ the $j$th coordinate vector in $\R^l$.
Then (\ref{EmanifU}) can be written as 
\[
\Pb(A_{all},\, \by(\bE,v_0) \in U,\, h(\by(\bE,v_0)) \cdot \bee_l=0)= p(M,U).
\]
As a consequence of the Taylor decomposition, for all $j\in\{0,\ldots,l-1\}$,
for any $\bu\in U$ and for all small $r$, 
\begin{equation}
\label{Etaylor}
h(\bu + r\bee_j)\cdot \bee_l= h(\bu) \cdot \bee_l+ r \, \rD h(\bu)\,\bee_j \cdot \bee_l + {\rm err}(\bu,j,r), 
\end{equation}
where for each $\bu \in U$, $\rD h(\bu)$ is the differential operator of $h$ at $\bu$, and
where the error term err$(\bu,j,r)=o(r)$ as $r\to 0$.
Since $h$ is a diffeomorphism, given any $\bu\in U$, there exists a 
$j\equiv j(\bu)\in \{0,\ldots,l-1\}$ 
such that
$\rD h(\bu) \,\bee_{j+1} \cdot \, \bee_l \neq 0$.
Therefore (\ref{EmanifU}) implies that, for some $j\in \{0,\ldots,l-1\}$,
\begin{equation}
\label{EmanifUj}
\Pb(A_{all},\, \by(\bE,v_0) \in M \cap U,\, \rD h(\by(\bE,v_0))\, \bee_{j+1} \cdot \bee_l> 0)\geq  \frac{p(M,U)}{2\,l},
\end{equation}
or 
\[
\Pb(A_{all},\, \by(\bE,v_0) \in M \cap U,\, \rD h(\bE(\bx,\iota_0))\, \bee_{j+1} \cdot \bee_l< 0)\geq  \frac{p(M,U)}{2\,l}.
\]
Without loss of generality, suppose (\ref{EmanifUj}) and choose $c,d\in (0,\infty)$, $c<d$,
and $\delta=\delta(c)>0$ such that 
\begin{eqnarray}
\nonumber
\Pb(A_{all},\, \by(\bE,v_0) \bck&\in&\bck M \cap U,\, \rD h(\by(\bE,v_0))\, \bee_{j+1} \cdot \bee_l\in (c,d),\\
\label{EmanifUjcd}
\bck\bck\bck& &\bck\bck\bck \sup_{r\in (0,\delta)}|{\rm err}(\by(\bE,v_0),j+1,r)|/r \leq c/2)
\geq  \frac{p(M,U)}{4\,l}.
\end{eqnarray}
Consider the modified processes $\eta_r^j(\bE)$, $r>0$,
corresponding to this $j$, and note that 
$\by(\eta_r^j(\bE),v_0)=\by(\bE,v_0)+r \,\bee_{j+1}$.
Now, due to (\ref{Etaylor}) and (\ref{EmanifUjcd}),
 one can choose a decreasing sequence $(r_m)_{m=1}^\iy$ of positive numbers converging to $0$, 
 so that the intervals defined by 
$J(r_m):=(c\, r_m/2, d\, r_m + c\, r_m/2)$, for each $m\geq 1$,
are mutually disjoint 
(i.e., $J(r_m)\cap J(r_{m'}) =\emptyset$  for $m < m'$) and such that
\[
\Pb(\eta_{r_m}^j(\bE)\in \bA_{all}, \, h(\by(\eta_{r_m}^j(\bE),v_0)) \cdot \bee_{j+1}\in J(r_m))\geq 
\frac{p(M,U)}{4\,l},
\]
hence
\[
\Pb(\bE\in \bA_{all}, \, h(\by(\bE,v_0)) \cdot \bee_{j+1}\in J(r_m))\geq 
e^{-r_m(W(0)+W(1))}
\frac{p(M,U)}{4\,l}.
\]
As in the linear case, one arrives to a contradiction.

\subsection{Odd versus even in the martingale technique}
The reason why the martingale technique fails on odd cycles is
similar: there is no non-trivial martingale that can be expressed as a
linear combination of the different  $W^*(X_n^{e_i})$,
$i=0,\ldots,l-1$, as in identity 
(\ref{eq:wstar}).
Indeed, let us fix a time $n\in\N$ and let, for all $i\in\Z/l\Z$, 
\begin{align*}
y_i&:= \Es(Y_{n+1}^+(i)-Y_{n}^+(i)|\F_n)=\Es(Y_{n+1}^-(i)-Y_{n}^-(i)|\F_n),\\
z_i&:=\Es(W^*(X_{n+1}^{e_i})-W^*(X_n^{e_i})|\F_n),
\end{align*}
$$
Y_n:=\left(\begin{array}{c}
y_0\\
\vdots\\
\
y_{l-1}
\end{array}
\right), ~
Z_n:=\left(
\begin{array}{c}
z_0\\
\vdots \\
z_{l-1}
\end{array}
\right).
$$
Then, for all $0\le i\le l-1$, 
$$z_i=y_i+y_{i+1},$$ 
since $W^*(X_{n+1}^{e_i})=Y_{n}^+(i)+Y_{n}^-(i+1)$. This
implies again that, almost surely,
$$Z_n=M^{(l)}Y_n, \ n\geq 1.$$
Suppose there is a fixed vector $\be\in \Real^l$ such that
the dot product $\be Y_n$ equals $0$, almost surely, for all $n$.
Since, at each time step $n\in\N$,  $Y_n$ has (almost surely) only one non-zero coordinate, namely,
$y_i>0$ for $i=I_n$ and $y_j=0$ for $j\neq I_n$, and since the walk visits each and every vertex
at least once with positive probability, we see that $\be$ is
necessarily the null-vector. 
As before, if $l$ is odd, $M^{(l)}$ is a regular matrix, and 
therefore  no martingale can be expressed 
as a non-trivial deterministic
linear combination of the different  $W^*(X_n^{e_i})$, $i=0,\ldots,l-1$.

However, we show in \cite{limtar} that, 
for all $i=0,\ldots,l-1$, if $t_n^i$ is the 
$n$-th return time to the vertex $i$, the process 
$W^*(X_{t_n^i}^{e_i})-W^*(X_{t_n^i}^{e_{i-1}})$ approximates a
martingale.
The accuracy of this approximation depends on the regularity of 
the weight function $W$, hence our argument requires 
technical assumptions on $W$. 
In particular, the main theorem in \cite{limtar} implies
(\ref{A})
for strongly {\ERRW s},
where $W$ is nondecreasing. 

Even though the time-lines technique is simpler in general, one cannot adapt it similarly, since it uses the independence of random variables and is therefore unstable with respect to small perturbation. 

\setlength{\textheight}{21cm}

\newpage
{\bf Acknowledgment.} P.T. would like to thank Christophe Sabot for an interesting discussion. We are very grateful to the referee for useful comments.

\end{document}